\def\versiondate{11 Dec.\ 2002}
\input math.macros
\input Ref.macros

\newif\ifsubmit

\checkdefinedreferencetrue
\continuousfigurenumberingtrue
\theoremcountingtrue
\sectionnumberstrue
\forwardreferencetrue
\citationgenerationtrue
\nobracketcittrue
\hyperstrue
\initialeqmacro

\bibsty{myapalike}

\def\F{{\cal F}}
\def\G{{\cal G}}

\def\verts{{\ss V}}

\def\edges{{\ss E}}
\def\ed#1{[#1]}  
\def\ded#1{\langle #1 \rangle}  

\def\t #1 {T_{#1}}   
\def\tu #1 {T^{#1}}    

\def\**{\sum_{n=1}^\infty {r_n \over |T_n|^2L_nL_{n-1}}}

\def\parent#1{\hat{#1}}  
\def\hm{{\ss Ray}}  
\def\hits{{\ss Perc}}  
\def\wsf{{\ss WSF}}
\def\fsf{{\ss FSF}}

\def\cmp{\tau}  
\def\flow{{\cal P}}  
\def\STAR{\bigstar}        
\def\CYCLE{\diamondsuit}
\def\rnchi#1{\eta^{#1}}
\def\TreePath#1|#2{[#1 \mid #2]}  

\def\BLPSusf{\ref b.BLPSusf/, hereinafter referred to as [BLPS01]%
  \global\def\BLPSusf{\htmllocref{BLPSusf}{[BLPS01]}}}

\ifsubmit
\pageno=0
\title{Change Intolerance in Spanning Forests}
\vfil
\centerline{Deborah Heicklen\footnote{${}^1$}{
Lockheed Martin M\&DS,
3200 Zanker Road M/S X75,
San Jose, CA 95134
\email{deborah.w.heicklen@lmco.com}
}
and
Russell Lyons\footnote{${}^2$}{
Department of Mathematics,
Indiana University,
Bloomington, IN 47405-5701
\email{rdlyons@indiana.edu}
and
School of Mathematics,
Georgia Institute of Technology,
Atlanta, GA 30332-0160
\email{rdlyons@math.gatech.edu}
}
}
\vfil
\noindent Running head: Change Intolerance in Spanning Forests
\vfil
\noindent Author to review proofs:
\begingroup\parskip=0pt\parindent=3truein\obeylines\frenchspacing
Russell Lyons
School of Mathematics
Georgia Institute of Technology
Atlanta, GA 30332-0160
\email{rdlyons@math.gatech.edu}
Tel. 404-894-5463
Fax 404-894-4409 
\endgroup
\vfil\eject
\fi

\def\firstheader{\eightpoint\ss\underbar{\raise2pt\line
    {{\it J. Theoret. Probab.} {\bf }
    \hfil Version of \versiondate}}}

\beginniceheadline

\vglue20pt

\title{Change Intolerance in Spanning Forests}

\author{Deborah Heicklen and Russell Lyons}

\bottom{Primary
60J45%
. Secondary
60J10, %
60J15, 
60J65, %
60K35 
.}
{Spanning forests.}
{Research partially supported by an NSF postdoctoral fellowship (Heicklen),
and NSF grant DMS-9802663 and the Miller Institute at the University of
California (Lyons).}

\abstract{
Call a percolation process on edges of a graph {\bf change intolerant} if
the status of each edge is almost surely determined by the status of the
other edges. We give necessary and sufficient conditions for change
intolerance of the wired spanning forest when the underlying graph is a
spherically symmetric tree. 
}

\bsection{Introduction}{s.intro}

An important pair of probability measures on graphs that has been
studied for the last 10 years is that of free and wired
uniform spanning forests. These measures have intimate connections to a
number of other probabilistic models, such as random walks, domino tilings,
random-cluster measures, and Brownian motion, as well as to some topics
outside of or only tangentially related to probability theory, such as
harmonic Dirichlet functions and $\ell^2$-Betti numbers.
Completely separately,
notions of insertion and deletion tolerance, also known as finite energy,
have been important for
decades in the study of percolation and other models in statistical
physics. 
Normally, one would expect spanning forest measures to be neither insertion
nor deletion tolerant, but we shall see that this is not always the case.

We now recall briefly the definitions of these spanning forest measures.
A comprehensive study of uniform spanning forests appears in \BLPSusf; a
survey appears in \ref b.Lyons:bird/.
The origin of these measures was the proof by
\ref b.Pemantle:ust/ of the following conjecture of Lyons:
If an infinite graph $G$ is exhausted by finite subgraphs $G_n$, then the
uniform distributions on the spanning trees of $G_n$ converge weakly to a
measure supported on spanning forests of $G$; by ``spanning forest'', we
mean a subgraph without cycles that contains every vertex.
We call this limit measure the {\bf free uniform spanning forest} measure
($\fsf$), since there is another natural construction where the exterior of
$G_n$ is identified to a single vertex (``wired'') before passing to the
limit.
This second construction, which we call the {\bf wired uniform spanning
forest} measure ($\wsf$), was implicit in Pemantle's paper and was made
explicit by \ref b.Hag:rcust/.
We shall be concerned primarily with $\wsf$ in this paper.

We next define the tolerance notions that we shall investigate.
Let  $G=(\verts,\edges)$ be a connected graph, where $\verts$ is the vertex set
and $\edges$ is the edge set. 
Let $\mu$ be a probability measure on $2^\edges := \{0,1\}^\edges$.
We think of an element $\omega \in \{0,1\}^\edges$ as the subset of $\edges$
where $\omega(e) = 1$, and refer to $\omega$ as a {\bf configuration}.
The {\bf status} of an edge $e$ in a configuration $\omega$ means simply
$\omega(e)$.
For a subset $K$ of $\edges$, let $\F(K)$ denote the $\mu$-completion of the
$\sigma$-field on
$2^\edges$ determined by the coordinates $e \in K$. We call $\mu$ {\bf
change intolerant} if the status of each edge is determined by the rest of
the configuration, i.e., if for all $e \in\edges$, we have $\F( \{ e \})
\subset \F(\edges \setminus \{ e \})$.
A useful equivalent definition is that for all $e \in\edges$, we have
$\mu_e\perp\mu_{\neg e}$, where $\mu_e$ denotes $\mu$ conditioned on
$\omega(e) = 1$ and restricted to $\F(\edges \setminus \{ e \})$, and
$\mu_{\neg e}$ denotes $\mu$ conditioned on $\omega(e) = 0$ and restricted
to $\F(\edges \setminus \{ e \})$.
One way that a measure $\mu$ may be far from change intolerant is to be {\bf
insertion tolerant}, meaning that for each $e$, we have $\mu[\omega(e) = 1
\mid \F(\edges \setminus \{ e \})] > 0$ a.s., or, equivalently, $\mu_{\neg
e} \ll \mu_e$.
Another way that a measure $\mu$ may be far from change intolerant is to be {\bf
deletion tolerant}, meaning that for each $e$, we have $\mu[\omega(e) = 0
\mid \F(\edges \setminus \{ e \})] > 0$ a.s., or, equivalently, $\mu_e \ll
\mu_{\neg e}$.
Since all trees are infinite $\wsf$-a.s., it is never the case that $\wsf$
is deletion tolerant. However, it is reasonable to ask whether $\wsf$ is
{\bf essentially deletion tolerant}, meaning that for each $e$, we have
$\mu[\omega(e) = 0 \mid \F(\edges \setminus \{ e \})] > 0$ a.s.\ on the
event that both endpoints of $e$ belong to infinite components
of $\omega \setminus \{ e \}$.

We shall show that there are trees $T$ for which $\wsf$ on $T$ is
change intolerant, as well as $T$ for which $\wsf$ is insertion
tolerant and essentially deletion tolerant. In fact, we shall give a
simple necessary and sufficient condition for a spherically symmetric tree
$T$ to have a change-intolerant $\wsf$.
Although it cannot happen for a spherically symmetric tree, there are
trees on which $\wsf$ is change intolerant at certain edges and insertion
tolerant at others: simply take two trees, one of which is change intolerant,
the other of which is insertion tolerant, and join them by identifying a
vertex of one with a vertex of the other.

There are many examples of graphs for which $\wsf$ is change intolerant.
Indeed, it is at first surprising that there should be any graphs for which
$\wsf$ is not change intolerant. It would be interesting to know the
situation on $T \times \Z$ when $T$ is a tree on which $\wsf$ is not change
intolerant.

\procl x.finite
Let $G$ be a finite graph. Then the uniform spanning tree is change
intolerant,
since the number of edges in a spanning tree is constant and equal to
$|\verts|-1$. 
\endprocl

For the next example we need some definitions. An infinite path in a
tree that starts at any vertex and does not backtrack is called a {\bf
ray}. Two rays are {\bf equivalent} if they have infinitely many
vertices in common. An equivalence class of rays is called an {\bf
end}.

\procl x.cayley
Let $G$ be the Cayley graph of a group which is not a finite extension
of $\Z$. Then by Theorem 10.1 of \BLPSusf, every component
tree in the $\wsf$ has exactly one end. Since adding an edge $e$ to a
configuration will either form a cycle or
connect two components, the resulting configuration would not be a forest
with trees having only one end each. Likewise, deleting an edge would
necessarily leave a finite component. Thus, $\wsf$ is change intolerant.
\endprocl

\procl x.twoZ3
Let $G$ be two copies of $\Z^3$ attached by an additional edge at
their respective origins. By Theorem 9.4 of \BLPSusf, the $\wsf$
has two components a.s., each with one end.
Hence adding an edge to a configuration would
either form a cycle or yield only one component, which cannot happen.
Similarly, deleting an edge would leave a finite component. Thus, $\wsf$ is
change intolerant.
\endprocl

More generally, $\wsf$ is change intolerant on any
graph for which the $\wsf$ has a finite number of components (since the
number of components is an a.s.\ constant by Theorem 9.4 of \BLPSusf)
or for which
each component has the same a.s.\ constant finite number of ends.
Also, recall that the $\wsf$ is a single tree a.s.\ on any graph for which
simple random walk is recurrent (\BLPSusf, Proposition 5.6).

In order to state our principal theorem, we need some notation.
Given a tree $T$ with root $o$ and given $k \in \N$, the {\bf level}
$\t k $ is the set of vertices of $T$ at distance $k$ from $o$. The tree 
is called {\bf spherically symmetric} if for each $k$, all vertices in
$\t k $ have the same number of children (the same degree).
It is well known that simple random walk on $T$ is transient iff $\sum_m
1/|\t m | < \infty$.

\procl t.main-default
Let $T$ be a spherically symmetric transient tree and
$L_n:= \sum_{m>n} 1/|\t m |$. 
Consider the series
$$
\sum_{n \ge 1}  {1 \over |\t n |^2 L_n L_{n-1}}
\,.
\label e.series
$$
\beginitems
\itemrm{(i)}
If this series diverges, then the $\wsf$ on $T$ is change
intolerant.
\itemrm{(ii)}
If this series converges, then the $\wsf$ on $T$ is insertion tolerant.
\itemrm{(iii)}
If this series converges and $T$ has bounded degree, then the $\wsf$ on $T$
is essentially deletion tolerant.
\enditems
\endprocl

Note that the series \ref e.series/ converges if ${|T_n|/n^\gamma}$ is bounded
above and below by positive constants for some $\gamma>1$, while the series
diverges if $\liminf_{n \to\infty} |T_{n+1}|/|T_n| \in (1, \infty)$.
Statement (iii) is not vacuous because when the series converges, a.s.\ all
components have more than one end by Corollary 11.4 of \BLPSusf, which is
reproduced as \ref c.sphsym/ below.

In the next section, we shall use the analysis of the $\wsf$ on trees given
by \BLPSusf\ to reduce \ref t.main-default/ to questions involving random
walk and percolation (\ref l.criteria/).
In \ref s.mart/, motivated by \ref b.LPP:concep/,
we introduce the martingales needed to
analyze these new questions and prove \ref t.main-default/.
(Actually, we extend \ref t.main-default/ to networks where the conductance
of an edge can depend on its distance to the root, as explained in \ref
s.prelim/.)
In the last section we give  an application to
determinantal probability measures and answer a question of \ref b.L:det/.

\bsection {Preliminary Reduction}{s.prelim}

If $x$ and $y$ are endpoints of an edge $\ed{x, y}$, we write $x \sim y$
and call $x, y$ {\bf neighbors}.
A {\bf network} is a pair $(G,C)$, where $G$ is
a connected graph with at least two vertices and $C$ is a function
from the unoriented edges of $G$ to the positive reals. The quantity
$C(e)$ is called the {\bf conductance} of $e$. The network is finite
if $G$ is finite. 
We assume that $\sum_{y \sim x} C(\ed{x, y})$ is finite for all
$x\in \verts$. 
The {\bf network random walk} on $(G, C)$ is the nearest-neighbor
random walk on $G$ with
transition probabilities proportional to the conductances.
The most natural network on $G$ is the default network
$(G,{\bf 1})$, for which the network random walk is simple random walk.
For each edge $e\in \edges$, the quantity $R(e) := 1/{C(e)}$ is the
{\bf resistance} of $e$. For general networks, we use a measure on spanning
trees adapted to the conductances. That is, for a finite network, choose
a spanning tree proportional to its weight, where the
{\bf weight} of a spanning tree $T$ is $\prod_{e\in T}C(e)$.
The proof in \ref b.Pemantle:ust/ of the existence of a limit of such measures
when an infinite graph is approximated by finite subgraphs, either
wired or not, extends to general networks. Explicit details are given
in \BLPSusf.

Given a tree $T$,
choose arbitrarily a vertex $o$ of $T$, which we call the {\bf root}.
If $x$ is a vertex of $T$ other than the root, write $\parent x$ for the
{\bf parent} of $x$, i.e., the next vertex after $x$ on the shortest path from
$x$ to the root.
We also call $x$ a {\bf child} of $\parent x$.
If $x$ and $y$ are two vertices, $x \wedge y$ denotes the most recent
common ancestor to $x$ and $y$.

Let $(T, C)$ be a transient network whose underlying graph is a tree.
For each vertex $x$, consider an independent network random walk $Z_x$
on $T$ starting at $x$ and stopped when it reaches $\parent x$, if ever.
(Note that $Z_o$ is never stopped.)
Write $\zeta_x$ for the set of edges crossed an odd number of times by
$Z_x$ and $\Xi$ for the set of vertices $x$ such that for all $y$,
if $x$ is an endpoint
of an edge in $\zeta_y$ and $\zeta_y$ is infinite, then $x = y$.
Then Section 11 of \BLPSusf\ shows that the law of $\bigcup_{x
\in\Xi} \zeta_x$ is $\wsf$.
In particular, if the deletion of an edge $e$ from $T$ leaves a recurrent
component, then $e$ belongs to the wired spanning forest a.s.
Therefore, we shall henceforth restrict consideration to networks that are
{\bf fully transient}, meaning that there are no edges whose deletion leaves a
recurrent component.

Note that $\zeta_o$ is a ray starting at $o$.
The law of $\zeta_o$ is denoted $\hm_o$.
The law of the connected component of $o$ in $\big\{ \ed{x, \parent x} \st x
\ne o,\, \zeta_x = \{ \ed{x, \parent x} \}\big\} $ is denoted $\hits_o$
(since the latter is an independent percolation on $T$).
Write $\hm_o \oplus \hits_o$ for the law of 
$\xi \cup \omega$ when $(\xi, \omega)$ has the law $\hm_o \otimes
\hits_o$.
By the description above, the component $\cmp(o)$ of $o$ in the $\wsf$ has
the same law as $\hm_o \oplus \hits_o$.

For neighbors $x, y$ in $T$, write $T_{x, y}$ for the component of $x$ in
$T \setminus \ed{x, y}$.
We take the root of $T_{x, y}$ to be $x$.
Let $\hm_{x, y}$ and $\hits_{x, y}$ denote $\hm_x$ and $\hits_x$, respectively,
on $T_{x, y}$.
Let $\cmp(x, y) := \cmp(x) \cap T_{x, y}$, where $\cmp(x)$ is the component
of $x$ in the $\wsf$ on all of $T$.
Let $\alpha_{x, y}$ be the probability that a network random walk on
$T$ started at $x$ is in $T_{x, y}$ from some time onwards.
The above alternative description of the $\wsf$ shows that
the $\wsf_{\neg \ed{x, y}}$-law of $\big(\cmp(x, y), \cmp(y, x)\big)$ is 
$$
(\hm_{x, y} \oplus \hits_{x, y}) \otimes (\hm_{y, x} \oplus \hits_{y, x})
\,,
\label e.law0
$$
while the $\wsf_{\ed{x, y}}$-law of $\big(\cmp(x, y), \cmp(y, x)\big)$ is 
the mixture 
$$
\alpha_{x, y} (\hm_{x, y} \oplus \hits_{x, y}) \otimes \hits_{y, x}
+
(1 - \alpha_{x, y}) \hits_{x, y} \otimes (\hm_{y, x} \oplus \hits_{y, x})
\,.
\label e.law1
$$
Note that $\alpha_{x, y} > 0$ since $T_{x, y}$ is transient and
$\alpha_{x, y} < 1$ since $T_{y, x}$ is transient.
Furthermore, the configuration of the $\wsf$ on the edges that do not touch
either $\cmp(x)$ or $\cmp(y)$ is independent of the status of $\ed{x, y}$
(given $\tau(x) \cup \tau(y)$).
We shall term the event $ \{ \omega \hbox{ is infinite} \}$ {\bf survival} by
analogy with branching processes.
This gives us the following criteria,
where $A_{x, y}:= \hits_{x, y}(\hbox{survival})$ and $\hits^*_{x, y}$ denotes
$\hits_{x, y}$ conditioned on survival when $A_{x, y} > 0$.

\procl l.criteria
Let $(T,C)$ be a fully transient network on a tree.
\beginitems
\itemrm{(i)}
$\wsf$ is change intolerant iff for all neighbors $x, y$,
$$
\hm_{x, y} \oplus \hits_{x, y} \perp \hits_{x, y}
\,.
$$
\itemrm{(ii)}
$\wsf$ is
insertion tolerant iff for all neighbors $x, y$,
$$
\hm_{x, y} \oplus \hits_{x, y} \ll \hits_{x, y}
\,.
$$
\itemrm{(iii)}
$\wsf$ is
essentially deletion tolerant iff for all neighbors $x, y$,
$$
A_{x, y} > 0 \quad\hbox{ and }\quad
\hits^*_{x, y} \ll \hm_{x, y} \oplus \hits_{x, y}
\,.
$$
\enditems
\endprocl

\proof
This is a straightforward comparison of \ref e.law0/ and \ref e.law1/, and so
we make only a few remarks on the proof. For
example, change intolerance is equivalent to the mutual singularity of \ref
e.law0/ and \ref e.law1/. Each term of the mixture \ref e.law1/ must be
singular to \ref e.law0/. Since $x$ and $y$ may be switched,
we may consider only the second term.
Since $\alpha_{x, y} < 1$, this singularity is 
$$
(\hm_{x, y} \oplus \hits_{x, y}) \otimes (\hm_{y, x} \oplus \hits_{y, x})
\perp
\hits_{x, y} \otimes (\hm_{y, x} \oplus \hits_{y, x})
\,,
$$
which is the same as 
$\hm_{x, y} \oplus \hits_{x, y} \perp \hits_{x, y}$.

In (iii), the word ``essentially" implies that we are concerned only with
the event that both $x$ and $y$ belong to infinite components when $\ed{x,
y}$ is removed from the wired spanning forest. This condition leads
therefore to consideration of $\hits^*_{x, y}$.
\Qed

Let $\flow(T)$ denote the set of unit flows on $T$ from $o$ to infinity,
i.e., the set of nonnegative functions $\theta$ on the vertices of $T$ such
that $\theta(o) = 1$ and for each vertex $x$, the sum of $\theta(y)$ over
all children $y$ of $x$ equals $\theta(x)$.
Let $h(x)$ be the probability that a network random walk starting at $x$
ever visits $o$.
Theorem 11.1 of \BLPSusf\ contains the following information:

\procl t.GenTrees 
Let $(T,C)$ be a transient network on a tree.
If for all $\theta \in \flow(T)$, the sum
$$
\sum_{x \neq o} \theta(x)^2 [h(x)^{-1}-h(\parent x )^{-1}]
\label e.ser1
$$
diverges, then all components
of the $\wsf$ on $T$ have one end a.s.; if this sum converges
for some $\theta \in \flow(T)$, then a.s.\ the $\wsf$ on $T$ has components
with more than one end.
\endprocl

A network on a spherically symmetric tree is itself called {\bf spherically
symmetric} if for all $k$, every edge connecting $\t k-1 $ with $\t k $ has
the same resistance, which we shall denote $r_k$.
Corollary 11.4 of \BLPSusf\ specializes \ref t.GenTrees/ to spherically
symmetric trees:

\procl c.sphsym 
Let $(T, C)$ be a spherically symmetric network on a tree.
Assume that the resulting network is transient,
i.e., $\sum_m  r_m/|\t m | < \infty$.
Denote $L_n:= \sum_{m>n} r_m/|\t m |$. If the sum
$$
\sum_{n \ge 1}  {r_n \over |\t n |^2 L_n L_{n-1}}
\label e.sphser1
$$
diverges, then all components of the $\wsf$ on $T$ have one end a.s.;
 if this series converges, then a.s.\ 
 all components of the $\wsf$ on $T$ have uncountably many ends. 
\endprocl

More specifically, the sum \ref e.sphser1/ is the minimum of the sum 
\ref e.ser1/ over $\flow(T)$; the minimum is achieved for the equally
splitting flow $\theta(x) := |\t |x| |^{-1}$, where $|x|$ denotes the distance
of $x$ to the root $o$.

\bsection{Martingales}{s.mart}

Write $I(x) := \hm_o\big( \{ \xi \st x \in \xi \} \big)$ for
the harmonic measure of the set of rays that pass through $x$.
For any subtree $t$ of $T$ containing $o$ and any $x \in t$,
write $\nu[t, x]$ for the probability that when $(\xi, \omega)$ has
the law $\hm_o \otimes \hits_o$, we have $x \in\xi$ and the first $|x|$ levels
of $\xi \cup\omega$ agree with those of $t$.
Since
$$
h(x) = \hits_o\big[ x \in \omega \big]
\,,
$$
we have
$$
\nu[t, x] = {I(x) \over h(x)} \hits_o([t]_{|x|})
\,,
$$
where $[t]_n$ denotes the set of subtrees of $T$ (rooted at $o$)
whose first $n$ levels agree with those of $t$.
Therefore, if $T$ is transient,
$$
(\hm_o \oplus \hits_o)([t]_n) = \sum_{x \in t_n} \nu[t, x]
= W_n(t) \hits_o([t]_n)
\,,
\label e.RN
$$
where 
$$
W_n(t) := \sum_{x \in t_n} {I(x) \over h(x)}
\,.
$$

Let ${\cal F}_n$ be the $\sigma$-algebra generated by the statuses of the
edges in $T$ whose endpoints are at distance at most $n$ from $o$.
Equation \ref e.RN/ says that
$W_n$ is the Radon-Nikod\'ym derivative of $\hm_o \oplus \hits_o$
restricted to $\F_n$ with respect to $\hits_o$
restricted to $\F_n$.
Therefore, the sequence $\Seq{W_n, \F_n}$ is a
martingale with respect to $\hits_o$.
Hence, $\Seq{W_n^2, \F_n}$ is a submartingale with respect to
$\hits_o$, which is the same thing as saying that $\Seq{W_n, \F_n}$ is a
submartingale with respect to $\hm_o \oplus \hits_o$.
Write $W(t) := \limsup_{n \to\infty} W_n(t)$.
We shall use the following general lemma (see
\ref b.Durrett:book/, Chapter 4, Theorem 3.3, p.~242, or \ref b.LP:book/).

\procl l.rndich Let $\kappa$ be a finite measure and $\lambda$ a probability
measure on a $\sigma$-field $\G$. Suppose that $\G_n$ are increasing
sub-$\sigma$-fields whose union generates $\G$ and that $(\kappa \restrict
\G_n)$ is absolutely continuous with respect to $(\lambda \restrict \G_n)$
with Radon-Nikod\'ym derivative $X_n$. Set $X := \limsup_{n\to\infty} X_n$.
Then
$$
\kappa \ll \lambda \iff X < \infty \quad \kappa\hbox{-a.e.}
$$
and
$$
\kappa \perp \lambda \iff X = \infty \quad \kappa\hbox{-a.e.}
$$
\endprocl

Write $A_o := \hits_o(\hbox{survival})$ and, when $A_o > 0$, let
$\hits^*_o$ be $\hits_o$ conditioned on survival.

\procl l.ourrn
If $(T, C)$ is fully transient, then
$$
\hm_o \oplus \hits_o \perp \hits_o  \iff W = \infty
\quad \hm_o \oplus \hits_o\hbox{-a.s.}
$$
and
$$
\hm_o \oplus \hits_o \ll \hits_o \iff W < \infty
\quad \hm_o \oplus \hits_o\hbox{-a.s.}
\label e.instol
$$
If also $A_o > 0$, then
$$
\hits^*_o \ll \hm_o \oplus \hits_o \iff W > 0 \quad
\hits^*_o\hbox{-a.s.}
\label e.deltol
$$
\endprocl

\proof
The first two statements follow from \ref l.rndich/ applied to $\kappa :=
\hm_o \oplus \hits_o$, $\lambda := \hits_o$, and $\G_n := \F_n$, since then
$X_n = W_n$.
The last statement follows from \ref l.rndich/ applied to $\kappa :=
\hits^*_o$, $\lambda := \hm_o \oplus \hits_o$, and $\G_n := \F_n$.
Indeed, by \ref e.RN/, we have
$$
W_n(t)^{-1} (\hm_o \oplus \hits_o)([t]_n) 
= \hits_o([t]_n)
$$
when $W_n(t) \ne 0$.
Our assumption of full transience implies that $W_n(t) \ne 0$ for all
infinite $t$.
Therefore, 
$$
\hits^*_o([t]_n)
= 
(A_o W_n(t))^{-1} (\hm_o \oplus \hits_o)([t]_n) 
\,.
$$
Thus, \ref l.rndich/ applies with $X_n = (A_o W_n)^{-1}$.
\Qed

For $x, y \in T_n$, write $\TreePath x|y$ for the set of edges between $x
\wedge y$ and $x$.
For $x \in T_n$, we have 
\begineqalno
(\hm_o \oplus \hits_o)[x \in t]
&=
\sum_{y \in T_n} \hm_o[y \in\xi] \hits_o\big[\TreePath x|y \subseteq
\omega\big]
\cr&=
\sum_{y \in T_n} I(y) {h(x) \over h(x \wedge y)}
\,.
\endeqalno
Therefore
\begineqalno
\int W_n \,d(\hm_o \oplus \hits_o)
&=
\int \sum_{x \in T_n} {I(x) \over h(x)} \II{x \in t}
\,d(\hm_o \oplus \hits_o)(t)
\cr&=
\sum_{x \in T_n} {I(x) \over h(x)} (\hm_o \oplus \hits_o)[x \in t]
\cr&=
\sum_{x, y \in T_n} {I(x) I(y) \over h(x \wedge y)}
\cr&=
\sum_{u \in T,\, |u| \le n} {1 \over h(u)} \sum_{x, y \in T_n, x \wedge y = u}
I(x) I(y)
\cr&=
\sum_{u \in T,\, |u| \le n} {1 \over h(u)} \left \{ I(u)^2 - \sum_{u =
\parent v, |v| \le n} I(v)^2 \right \}
\cr&=
1 + \sum_{u \in T,\, 0 < |u| \le n} I(u)^2 \left \{{1 \over h(u)} - {1 \over
h(\parent u)} \right \}
\,.
\label e.Wnbound
\endeqalno
Note that these are the same summands that appear in \ref e.ser1/. 
Since $I(\cdot)$ is the equally splitting flow if $T$ is spherically
symmetric, it follows that
for a spherically symmetric tree, \ref e.Wnbound/ is bounded (in $n$)
iff \ref e.sphser1/ converges.

\procl t.main
Let $(T, C)$ be a spherically symmetric transient network on a tree.
\beginitems
\itemrm{(i)}
If the series \ref e.sphser1/ diverges, then the $\wsf$ on $T$ is change
intolerant.
\itemrm{(ii)}
If this series converges, then the $\wsf$ on $T$ is insertion tolerant.
\itemrm{(iii)}
If this series converges and $T$ has bounded degree, then the $\wsf$ on $T$
is essentially deletion tolerant.
\enditems
\endprocl

\proof
Note that the spherical symmetry implies that $(T, C)$ is fully transient.

(i) In this case, \ref c.sphsym/ shows that each component of the $\wsf$
has only one end a.s., whence no edge can be either inserted nor deleted.

(ii) When \ref e.sphser1/ converges, \ref e.Wnbound/ shows that $\Seq{W_n}$
is bounded in expectation with respect to $\hm_o \oplus \hits_o$.
Since $\Seq{W_n}$ is a submartingale, it follows that $W < \infty$ a.s.\
with respect to $\hm_o \oplus \hits_o$. 
In view of \ref e.instol/, we obtain $\hm_o \oplus \hits_o \ll \hits_o$. 
In particular, $A_o > 0$, which we shall use in proving (iii).
Now for any neighbors $x, y \in T$, the subtree $T_{x, y}$ is composed of a
finite collection of spherically symmetric trees attached at the leaves of a
finite tree.
Thus, a similar argument shows that for all neighbors $x, y$,
we have $\hm_{x, y} \oplus \hits_{x, y} \ll \hits_{x, y}$.
Therefore, the $\wsf$ on $T$ is insertion tolerant by \ref l.criteria/.

(iii) Another way to regard $\hits_o$ is as a branching process in a
varying environment (BPVE). 
With this view, $\Seq{W_n, \F_n}$ is the usual martingale in the theory of
branching processes.
Theorem 4.14 of \ref b.Lyons:rwcpt/ shows that when the offspring
distribution of a BPVE is uniformly bounded, then $W > 0$ a.s.\ given the
event of survival. Since $A_o > 0$ (proved for part (ii)),
this proves that when $T$ has bounded degree, $W >
0$ a.s.\ with respect to $\hits_o^*$. Because of \ref e.deltol/, we obtain
$\hits^*_o \ll \hm_o \oplus \hits_o$. As for part (ii), a similar argument
shows that for all neighbors $x, y$, we have $\hits^*_{x, y} \ll \hm_{x, y}
\oplus \hits_{x, y}$.
Therefore, the $\wsf$ on $T$ is essentially deletion tolerant by \ref
l.criteria/.
\Qed



\bsection{Singularity of Determinantal Probabilities}{s.det}

Let $(G, C)$ be a finite or infinite network. For this section, we choose
an orientation for each edge.
Identify each $e \in \edges$ with the corresponding unit vector $\I e$ in
$\ell^2(\edges)$.
Given two neighbors $x, y$, let
$$
\rnchi{\ded{x, y}}:= \cases{\ded{x, y} &if $\ded{x, y} \in \edges$\cr
                     - \ded{x, y} &if $\ded{y, x} \in \edges$.\cr
                     }
$$
Let $\STAR$ denote the closure in $\ell^2(\edges)$
of the linear span of the 
{\bf stars} $\sum_{y \sim x} \sqrt{C(\ed{x, y})} \rnchi{\ded{x, y}}$
($x \in \verts(G)$).
For a cycle of vertices $x_0, x_1, \ldots, x_n = x_0$, the function
$$
\sum_{i=0}^{n-1} \rnchi{\ded{x_i, x_{i+1}}}/\sqrt{C(\ed{x_i, x_{i+1}})}
$$
is called a {\bf cycle}.
Let $\CYCLE$ be the closure of the linear span of the cycles.
Since each star and cycle are orthogonal to each other, we have
$\STAR \perp \CYCLE$. 

Given any subspace $H\subseteq\ell^2(\edges)$, let $P_H$
denote the orthogonal projection of $\ell^2(\edges)$ onto
$H$, and let $P^\perp_{H}$ denote the orthogonal
projection onto the orthogonal complement $H^\perp$ of $H$.
The following result of \BLPSusf\ (Theorem 7.8 in an isomorphic form)
extends the Transfer Current Theorem of \ref b.BurPem/:

\procl t.transfer
Given any network $G$ and
any distinct edges $e_1,\ldots,e_k\in G$, we have
$$
\fsf[\omega(e_1) = 1, \ldots, \omega(e_k) = 1] =
\det [(P^\perp_{\CYCLE}{e_i}, {e_j})]_{1\le i,j\le k}
$$
and
$$
\wsf[\omega(e_1) = 1, \ldots, \omega(e_k) = 1] =
\det [(P_\STAR {e_i}, {e_j})]_{1\le i,j\le k}
\,.
$$
\endprocl

Clearly, these formulas characterize $\fsf$ and $\wsf$. In particular, as
observed in \BLPSusf, $\STAR \subseteq \CYCLE^\perp$, with equality iff $\wsf
= \fsf$. Question 15.11 of \BLPSusf\ asks whether $\wsf \perp \fsf$ when
the two measures are not equal; some cases where this is known to be true
are stated there. This question remains open, but it suggested a more
general possibility to \ref b.L:det/, which we may now show is false.

First, we give the more general context in which the question arose.
Given any countable set $E$, identify each $e \in E$ with the corresponding
unit vector $\I e$ in $\ell^2(E)$.
Given any closed subspace $H \subset \ell^2(E)$,
there is a unique probability measure $\P^H$ on $2^E$ defined by 
$$
\P^H[\omega(e_1) = 1, \ldots, \omega(e_k) = 1]
=
\det [(P_H e_i, e_j)]_{1\le i,j\le k}
$$
for any set of distinct $e_1,\ldots,e_k\in E$; see \ref b.L:det/ and \ref
b.DVJ/, Exercises 5.4.7--5.4.8.
In case $H$ is finite dimensional, then $\P^H$ is concentrated on subsets
of $E$ of cardinality equal to the dimension of $H$.

This suggests that in general, if $H_1 \subset H_2 \subset \ell^2(E)$
and $H_1 \ne H_2$, then $\P^{H_1} \perp \P^{H_2}$, the question asked in
\ref b.L:det/. But this is false. To see how this follows from \ref
t.main-default/, we must consider the effect of conditioning on the
measure $\wsf$ and its representation via determinants.
This is done partly in \BLPSusf\ and fully in \ref b.L:det/. The result is
that if we identify $\ell^2(\edges \setminus \{ e \} )$ with $(\R e)^\perp
\subset \ell^2(\edges)$, then
$\wsf_e = \P^{H_1}$ and $\wsf_{\neg e} = \P^{H_2}$, where
$$
H_1 := \STAR \cap (\R e)^\perp
\quad\hbox{and}\quad
H_2 := (\STAR + \R e) \cap (\R e)^\perp
\,.
$$
Thus, $H_1 \subseteq H_2$; furthermore, $H_1 \ne H_2$ as long as
$e \notin \STAR$, i.e., $\wsf[\omega(e) = 1] < 1$.
This condition holds on a tree $T$ when $e = \ed{x, y}$ and both $T_{x, y}$
and $T_{y, x}$ are transient.
Yet $\P^{H_1} \perp \P^{H_2}$ does not hold
when $\wsf$ is insertion tolerant (at $e$), as it may be.

\medbreak
\noindent {\bf Acknowledgement.}\enspace We are grateful to the referee for
a careful reading and useful suggestions.

\bibfile{\jobname}
\def\noop#1{\relax}
\input \jobname.bbl

\filbreak
\begingroup
\eightpoint\sc
\parindent=0pt\baselineskip=10pt

Lockheed Martin M\&DS,
3200 Zanker Road M/S X75,
San Jose, CA 95134
\email{deborah.w.heicklen@lmco.com}

Department of Mathematics,
Indiana University,
Bloomington, IN 47405-5701
\emailwww{rdlyons@indiana.edu}
{http://php.indiana.edu/\string~rdlyons/}

and

School of Mathematics,
Georgia Institute of Technology,
Atlanta, GA 30332-0160
\email{rdlyons@math.gatech.edu}

%

\endgroup


\def\temp{\let\linkit=\linkyear \apaliketrue}
\temp
\ifcitationgeneration\immediate\write\labelfile{\sanitize\temp}\fi
\def\startreferences{
 \vskip0pt plus.3\vsize \penalty -150 \vskip0pt
 plus-.3\vsize \bigskip\bigskip \vskip \parskip
 \begingroup\baselineskip=12pt\frenchspacing
 \bibliographytitle
 \vskip12pt\parindent=0pt
 \def\and{{\rm and}}
 \def\em{\it}
 \def\newblock{\hskip .11em plus.33em minus.07em}
 \def\bibauthor##1{{\sc ##1}}
 \def\bibitem[##1]##2
 {\htmlanchor{##2}{}\RefLabel{##2}[##1]\hangindent=.8cm\hangafter=1}
 }
\def\endreferences{\bigskip\bigskip\endgroup}
\ifundefined{bibstylemodification}\relax\else\bibstylemodification\fi
\startreferences

\bibitem[Benjamini, Lyons, Peres, and Schramm (2001)]{BLPSusf}
\bibauthor{Benjamini, I., Lyons, R., Peres, Y., \and{} Schramm, O.} (2001).
\newblock Uniform spanning forests.
\newblock {\em Ann. Probab.} {\bf 29}, 1--65.

\bibitem[Burton and Pemantle (1993)]{BurPem}
\bibauthor{Burton, R.M. \and{} Pemantle, R.} (1993).
\newblock Local characteristics, entropy and limit theorems for spanning trees
  and domino tilings via transfer-impedances.
\newblock {\em Ann. Probab.} {\bf 21}, 1329--1371.

\bibitem[Daley and Vere-Jones (1988)]{DVJ}
\bibauthor{Daley, D.J. \and{} Vere-Jones, D.} (1988).
\newblock {\em An Introduction to the Theory of Point Processes}.
\newblock Springer-Verlag, New York.

\bibitem[Durrett (1996)]{Durrett:book}
\bibauthor{Durrett, R.} (1996).
\newblock {\em Probability: Theory and Examples}.
\newblock Duxbury Press, Belmont, CA, second edition.

\bibitem[H{\"a}ggstr{\"o}m (1995)]{Hag:rcust}
\bibauthor{H{\"a}ggstr{\"o}m, O.} (1995).
\newblock Random-cluster measures and uniform spanning trees.
\newblock {\em Stochastic Process. Appl.} {\bf 59}, 267--275.

\bibitem[Lyons (1992)]{Lyons:rwcpt}
\bibauthor{Lyons, R.} (1992).
\newblock Random walks, capacity and percolation on trees.
\newblock {\em Ann. Probab.} {\bf 20}, 2043--2088.

\bibitem[Lyons (1998)]{Lyons:bird}
\bibauthor{Lyons, R.} (1998).
\newblock A bird's-eye view of uniform spanning trees and forests.
\newblock In Aldous, D. \and{} Propp, J., editors, {\em Microsurveys in
  Discrete Probability}, volume 41 of {\em DIMACS Series in Discrete
  Mathematics and Theoretical Computer Science}, pages 135--162. Amer. Math.
  Soc., Providence, RI.
\newblock Papers from the workshop held as part of the Dimacs Special Year on
  Discrete Probability in Princeton, NJ, June 2--6, 1997.

\bibitem[Lyons (2002)]{L:det}
\bibauthor{Lyons, R.} (2002).
\newblock Determinantal probability measures.
\newblock Preprint.

\bibitem[Lyons, Pemantle, and Peres (1995)]{LPP:concep}
\bibauthor{Lyons, R., Pemantle, R., \and{} Peres, Y.} (1995).
\newblock Conceptual proofs of ${L}\log {L}$ criteria for mean behavior of
  branching processes.
\newblock {\em Ann. Probab.} {\bf 23}, 1125--1138.

\bibitem[Lyons with Peres (2003)]{LP:book}
\bibauthor{Lyons, R. {\rm with} Peres, Y.} (2003).
\newblock {\em Probability on Trees and Networks}.
\newblock Cambridge University Press.
\newblock In preparation. Current
  version available at \hfil\break
  \htmlref{http://php.indiana.edu/\string~rdlyons/}{{\tt
  http://php.indiana.edu/\string~rdlyons/}}.

\bibitem[Pemantle (1991)]{Pemantle:ust}
\bibauthor{Pemantle, R.} (1991).
\newblock Choosing a spanning tree for the integer lattice uniformly.
\newblock {\em Ann. Probab.} {\bf 19}, 1559--1574.

\endreferences
\bye